\documentclass[10pt]{amsart}
\usepackage{amsfonts}
\usepackage{amssymb}
\usepackage{amsmath}
\usepackage[arrow,matrix,curve]{xy}

\usepackage[active]{srcltx}

\usepackage[normalem]{ulem}
\usepackage{color}

\newcommand{\C}{{\mathbb C}}

\newcommand{\N}{{\mathbb N}}

\newtheorem{theorem}{Theorem}[section]
\newtheorem{lemma}[theorem]{Lemma}

\newtheorem{corollary}[theorem]{Corollary}

\newtheorem{definition}[theorem]{Definition}

\def\cal{\mathcal}

\newcommand{\calh}[0]{{\cal H}}

\newcommand{\calb}[0]{{\cal B}}
\newcommand{\cala}[0]{{\cal A}}

\newcommand{\calk}[0]{{\cal K}}

\newcommand{\cald}[0]{{\cal D}}

\newcommand{\cals}[0]{{\cal S}}

\begin{document}

\title{The generators and relations picture of $KK$-theory}
\author[B. Burgstaller]{B. Burgstaller}
\email{bernhardburgstaller@yahoo.de}
\address{}
\keywords{$KK$-theory, picture, generators}

\begin{abstract}
This is half an overview article since what we describe here is essentially known.
We describe $KK$-theory by generators and relations in a formal sum of 
formal products of $*$-homomorphisms and some synthetical morphisms.
What comes out is a category.
The Kasparov product is then just the composition of morphisms.
Our description may be interesting to anyone who wants a quick and elementary definition of $KK$-theory.
This description could also be used 
for other categories of algebras than $C^*$-algebras endowed with group actions, for example, $C^*$-algebras equipped with an action by a semigroup, a category et cetera.
\end{abstract}


\maketitle

\section{Introduction}

In 1980 G.G. Kasparov introduced $KK$-theory in his influential paper \cite{kasparov1980rus,kasparov1981}. Not less influential was the progress in his paper \cite{kasparov1988}, where
beside a wealth of new ideas
the difficulties in technicalities in $KK$-theory were reduced by the incorporation of results published in the meanwhile by A. Connes, G. Skandalis, and N. Higson.
J. Cuntz found out some universal aspect of $KK$-theory \cite{cuntz1984}, and these findings were elaborated further and brought to its final form by N. Higson
\cite{higson}.
Based on his findings, J. Cuntz found another 
picture of $KK$-theory \cite{cuntz1987}, which tends to be somewhat easier 
than Kasparov's technical approach.
Up to now, $KK$-theory showed up an immense impact in operator theory, $K$-theory, geometry, analysis and related topics like dynamical systems. 

In this note we describe another picture of $KK$-theory, which is based on the universal description 
by Higson \cite{higson}, and essentially quite clear.
Even if it seems evident, it has not yet been formulated in the literature to our best knowledge. It is based on generators and relations subject to relations dictated by
the universal properties of $KK$-theory. The generators are the $C^*$-homomorphisms together with certain synthetical inverses. 
We consider the formal sums of their formal products,
and in this free construction we introduce the relations to get a theory called $GK$-theory. 
It has the same universal property as $KK$-theory. 
For separable $C^*$-algebras 
$GK$-theory and $KK$-theory evidently coincide up to isomorphism.

The advantage of 
this approach is that it is quite elementary, and the interested reader needs only basic knowledge in $C^*$-theory and category theory for reading it.
In this way it may serve as a fast and easy study of the definitions of $KK$-theory. The technical Kasparov product is automatically included and need not further be studied, since it is implicitly given
by the definition of composition of morphisms in a category.
The reader who expects compressed mathematics will recognize that
we keep our exposition easy and light.
This is no accident as we want to address a possibly large audience.
In particular, we have also physicists in mind who want to have a quick but complete definition of $KK$-theory.

The reader who wants to understand the complete current literature in $KK$-theory will not come around the $KK$-theory picture by Kasparov, since it is mostly used.
Still, the generators and relations picture may also be interesting for those who are already familiar with $KK$-theory.

Another benefit of the generators and relations approach is that it works also for other categories of $C^*$-algebras, for example
the category of $C^*$-algebras endowed with an action by a semigroup, a category or whatsoever. The homomorphisms should then be understood to be
equivariant in the respective category. Also completely different algebras than $C^*$-algebras may be considered. The algebra $\calk$ of compact operators
has then probably to be substituted by another stabilizing algebra, say the closure $\overline{M_\infty}$ of the matrices $M_\infty$ under some topology.
The algebra of smooth compact operators in the approach by Cuntz \cite{cuntz1997} may come into mind.

Section \ref{section_def_gk} introduces the generators and relations picture of $KK$-theory, and requires only basic knowledge in $C^*$-theory and category theory.
In Section \ref{section_universal_property} the universal property of $GK$-theory is formulated.
Section \ref{section_homsets} is an addendum where morphism classes are turned to morphism sets. This works as soon as we introduce a bound for the cardinality of the 
$C^*$-algebras under consideration.

Section \ref{section_kk} shows that $GK$-theory and $KK$-theory coincide. It uses deep results of other authors
(as mentioned by J. Cuntz \cite{cuntz1984} and N. Higson \cite{higson}, and by K. Thomsen \cite{thomsen} and R. Meyer \cite{meyer2000}
for including group actions).
This section requires that the reader has some familiarity with $KK$-theory,
but is otherwise also 
easy.



\section{The definition of $GK$-theory}  

\label{section_def_gk}

\subsection{Some notations}

Let $M$ be a locally compact second countable group.
We shall consider the category $C^*$ with object class $Obj(C^*)$ consisting of all $C^*$-algebras $A$ endowed with an action by $M$.
The morphism set $C^*(A,B)$ from object $A$ to object $B$ is defined to be the set of all $M$-equivariant $*$-homomorphisms $f: A \rightarrow B$ from  $A$ to $B$.
Their collection forms the morphism class $Mor(C^*)$.
The letter $1_A$ denotes the identity morphism in $C^*(A,A)$.

The $C^*$-algebra of compact operators on a separable Hilbert space 
is denoted by $\calk$.
A {corner embedding} is a morphism $f: A \rightarrow A \otimes \calk$ in $C^*$ of the form $f(a)=a \otimes e$ for some one-dimensional $M$-invariant projection $e\in \calk$.
The algebra $A \otimes \calk$ may be endowed with any $M$-action, and need not be diagonal.  

\subsection{Motivation}
The reader not familiar with $KK$-theory may wish to skip this subsection in a first reading.

$KK$-theory means the category $KK$ with object class $Obj(C^*)$ and morphism set from $A$ to $B$ to be Kasparov's $KK$-theory group $KK(A,B)$.
Composition of morphisms is given by the Kasparov product.

We are going to describe a category $GK$, which is based on the universal description of $KK$-theory by Higson in \cite[Thm. 4.5]{higson},
namely that a certain functor $C^* \rightarrow KK$ ($KK$-functor) from $C^*$-theory to $KK$-theory is a universal functor into an additive category
which is homotopy invariant, stable and split-exact. For the definition of these properties see \cite[(i)-(iii) on page 269]{higson}.

To recall them we note that a functor $G:C^* \rightarrow E$ into an additive category $E$ (see \cite{maclane})
is said to be
\begin{itemize}
\item
stable, if every corner embedding $f: A \rightarrow A \otimes \calk$ in $C^*$ induces an isomorphism $G(f)$,

\item
homotopy invariant, if all homotopic morphisms $f_0,f_1:A \rightarrow B$ in $C^*$ induce an identity $G(f_0)=G(f_1)$, and

\item
split-exact, if for every split exact sequence (\ref{splitexact}) in $C^*$ the morphism $\eta$ defined in (\ref{defxi}) and entered in the diagram
(\ref{diagramG}) is invertible.

\end{itemize}

These properties say that the $KK$-theory functor ignores non-commutativity and works out the commutative structure of $C^*$.
Stability of the functor annihilates with $\calk$ the simplest purely non-commutative $C^*$-algebra (beside matrix algebras). 
A split exact sequence reflects roughly a direct sum of $C^*$-algebras, which is directly transported into the additive structure
in $KK$-theory by the split-exactness of the $KK$-theory functor. Only homotopy invariance tends to simplify also a commutative
context.

\subsection{Introducing new invertible morphisms}

\label{uu}

Our starting point is the category $C^*$ itself.  
In our construction we shall need to turn certain morphisms into invertible morphisms.
To this end we shall enrich the alphabet of homomorphisms $Mor(C^*)$ with a new collection of morphisms which will
be later defined to be inverses.

For all objects $A,B$ in $C^*$ define $\Theta(A,B)$ to be the set $C^*(A,B)$.
We are now enlarging these sets $\Theta(A,B)$ as follows.

\begin{itemize}

\item 
To all objects $A$ and $B$ in $C^*$ and every corner embedding (homomorphism)
$f:A \rightarrow A \otimes \calk$ in $C^*$ we add a new letter $f^{-1}$ to $\Theta(A \otimes \calk, A)$,
which will be later 
the inverse for $f$.

\item
For every split exact sequence
\begin{equation}  \label{splitexact}
\cals:\xymatrix{0 \ar[r] & A \ar[r]^{f} & D \ar@<.5ex>[r]^{g} & B \ar[r] \ar@<.5ex>[l]^{s} & 0  }
\end{equation}
in $C^*$
we add a new letter $\vartheta_\cals$ to $\Theta(D, A \oplus B)$.

\end{itemize}

After these enlargements, each $\Theta(A,B)$ is still a set.

\subsection{Introducing composition of morphisms}

In order to be able to compose our new morphisms $\Theta(A,B)$ we need a product
(which will later be the Kasparov product).
We shall choose the free product.

For all objects $A$ and $B$ in $C^*$,
let $\Lambda(A,B)$ be the class consisting of all finite sequences (free words)
$$f_1  f_2 f_3 \ldots  f_n$$
for which there exist 
objects $A_1,A_2,\ldots, A_{n+1}$ such that $A=A_1$, $B=A_{n+1}$ and
$f_i \in \Theta(A_i,A_{i+1})$.
We visualize 
$f_1 f_2 \ldots f_n$ as a path of morphisms like this:
\begin{equation}   \label{path}
\xymatrix{A_1 \ar[r]^{f_1} & A_2 \ar[r]^{f_2} & A_3 \ar[r]^{f_3} &\ldots \ar[r]^{f_n} & A_{n+1} . }
\end{equation}
Notice the reversed order, that $f_1 f_2 \ldots f_n$ will stand for the composition $f_n \circ f_{n-1} \circ
\ldots f_2 \circ f_1$ in standard notations for composition of homomorphisms.


The product
$$\Lambda(A,B) \times \Lambda(B,C) \rightarrow \Lambda(A,C)$$
is given formally by concatenation:
$$f_1 f_2 \ldots f_n \times g_1 g_2 \ldots g_m = f_1 f_2 \ldots f_n g_1 g_2 \ldots g_m$$
for $f_1 f_2 \ldots f_n \in \Lambda(A,B)$ and $g_1 g_2 \ldots g_m \in \Lambda(B,C)$.

In other words, $\Lambda$ may be visualized as the directed graph with vertices $Obj(C^*)$ and edges $\Theta$. 
The product in this graph is given by concatenation of paths.

\subsection{Introducing addition}

To obtain an additive category, we need to be able to add and subtract morphisms in $\Lambda(A,B)$. To this end,
for all objects $A$ and $B$ in $C^*$ we introduce
the class $\Gamma(A,B)$ consisting of all formal sums
\begin{equation}  \label{formalsum}
\pm f_1 \pm f_2 \pm \ldots \pm f_n
\end{equation}
where all $f_i \in \Lambda(A,B)$ and each $\pm$-sign stands here for either a $+$ or a $-$.
Addition and subtraction in $\Gamma(A,B)$ is given formally by concatenation of two expressions.
The multiplication defined in $\Lambda$ will be extended to a multiplication in $\Gamma$
by the distributive law. That is, for all objects $A,B,C$ in $Obj(C^*)$ we define multiplication
$$\Gamma(A,B) \times \Gamma(B,C) \rightarrow \Gamma(A,C)$$
as
\begin{eqnarray}   \label{defproduct}
&&(\pm f_1 \pm f_2 \pm \ldots \pm f_n) \times (\pm g_1 \pm g_2 \pm \ldots \pm g_m) \\
&=& \pm f_1 g_1 \pm \ldots \pm f_n g_1 \pm f_1 g_2 \pm \ldots \pm f_n g_2 \pm \ldots \ldots \pm f_1 g_m \pm  \nonumber
\ldots \pm f_n g_m.
\end{eqnarray}
The $\pm$-signs in the product are choosen as usual, that is, for example $(+f_1)\times (-g_1)=- f_1 g_1$,
$(-f_1)\times (-g_1)=+ f_1 g_1$ et cetera.


\subsection{Introducing equivalences}

\label{introducing_equivalences}

The substance of our construction is finished with the classes $\Gamma(A,B)$. Now we are going
to divide out relations in $\Gamma$ to turn it into our desired $KK$-category.

For all objects $A$ and $B$ in $C^*$
we say that two elements $f$ and $g$ in $\Gamma(A,B)$ are {\em equivalent} if there is a finite
sequence $f=f_1,f_2,\ldots, f_n =g$ in $\Gamma(A,B)$ such that two neighboring elements $f_i,f_{i+1}$
distinguish from each other by an {\em elementary equivalence} (or modification).

The first elementary equivalences that we shall introduce are those that turn the formal sums
in $\Gamma(A,B)$ into a real sum. 
That is, we allow as elementary equivalences in $\Gamma(A,B)$

\begin{itemize}
\item
the permutation of two neighboring summands (together with their signs) in (\ref{formalsum}),

\item
the cancelation of two neighboring elements within the expression (\ref{formalsum}), that is, $f-f \equiv 0$, where $0$ denotes the zero homomorphism
in $C^*(A,B)$,

\item
and the zero element relation $f+0 \equiv f$.

\end{itemize}

%
%
Each elementary equivalence
\begin{equation}   \label{elementaryequivalence}
f_1 \equiv f_2
\end{equation}
in $\Gamma(A,B)$
we have just introduced, like $f+g\equiv g+f$ for permutation of summands, and that we shall introduce
is understood that it can appear also in any composed expression of the form
\begin{equation}  \label{composedexpression}
x + y f_1 z \equiv x + y f_2 z
\end{equation}
for some $x \in \Gamma(D,E)$, $y \in \Gamma(D,A)$ and $z \in \Gamma(B,E)$
to form another 
elementary equivalence 
of the elementary equivalence (\ref{elementaryequivalence}).
Each single $x,y$ or $z$ can here also  not appear.


\subsection{Turning into an additive category}

To turn $\Gamma$ into an additive category, we need to introduce further equivalences.

In an Ab-category there exists the notion of a biproduct, see \cite[VIII.2.Def.]{maclane},
and the existence of a biproduct is equivalent to the existence of a product, or a coproduct,
which then are the same object, see \cite[VIII.2.Thm. 2]{maclane}.

We thus use the notion of the biproduct as an elementary equivalence as follows:

\begin{itemize}

\item
For all objects $A$ and $B$ in $C^*$ we have a diagram
\begin{equation}   \label{splitexactAB}
\xymatrix{A \ar@<.5ex>[r]^{i_A} & A \oplus B \ar@<-.5ex>[r]_{p_A} \ar@<.5ex>[l]^{p_A} & B \ar@<-.5ex>[l]_{i_B},  }
\end{equation}
where $A \oplus B$ denotes the $C^*$-algebraic direct sum of $A$ and $B$, $i_A,i_B,p_A,p_B$
the canonical injections and projections in $C^*$, and we introduce
the relation (elementary equivalence)
\begin{equation}   \label{biproduct}
p_A i_A +  p_B i_B \equiv 1_{A \oplus B}
\end{equation}
in $\Gamma(A\oplus B, A \oplus B)$,
where $1_{A \oplus B}$ denotes the identical morphism in $C^*$.
To get all elementary equivalences, it is understood that this relation can appear
also in a composed expression like in (\ref{composedexpression}).

\end{itemize}
(Notice that $p_A i_A$ will mean the composition $i_A \circ p_A$ of maps in standard notation.
We use however the reversed order in $\Gamma$.)

\subsection{Respecting homomorphisms}

\label{respecting_homomorphisms}

In order to get finally a functor from $C^*$ to our aimed category we need a further elementary
equivalence in $\Gamma$.

\begin{itemize}

\item
For all objects $A,B,C$ in $C^*$ and all morphisms $f:A \rightarrow B$ and $g:B \rightarrow C$
in $C^*$
we introduce the elementary equivalence
\begin{equation}    \label{equivalencemorphism}
 f g \equiv (g \circ f)
\end{equation}
in $\Gamma(A,C)$,
where $g \circ f$ denotes the composition of homomorphisms in $C^*$.
Again, it is understood that this equivalence may also appear in a composed expression
as in (\ref{composedexpression}).




\end{itemize}

\subsection{The unit elements}

\label{the_unit_element}

To obtain a category, we need 
identity morphisms in $\Gamma$.
To this end we need to take care
of the letters in $\Theta$ not appearing in $C^*$.

\begin{itemize}

\item
For every synthetical letter $f \in \Theta(A,B) \backslash C^*(A,B)$ introduced in \ref{uu} we introduce the elementary
equivalences
$$ 1_A f \equiv f$$
and
$$ f 1_B \equiv f$$
in $\Gamma(A,B)$.

\end{itemize}


\subsection{Respecting homotopy}

\label{respecting_homotopy}

We now turn to the defining relations of $KK$-theory. One is homotopy equivalence.

\begin{itemize}

\item
If $f:A \rightarrow B[0,1]$ is a homotopy between $f_0:A \rightarrow B$ and $f_1:A \rightarrow B$
in $C^*$
then we introduce an elementary equivalence
$$f_0 \equiv f_1$$
in $\Gamma(A,B)$.


\end{itemize}

\subsection{Respecting stability}

\label{respecting_stability}

Another characteristic of $KK$-theory is stability.

\begin{itemize}

\item
For all objects $A$ and $B$ in $C^*$ and every corner embedding $f:A \rightarrow A \otimes \calk$
in $C^*$ and $f^{-1}$ as
in \ref{uu} we introduce the elementary equivalences
$$f f^{-1} \equiv 1_A$$ 
in $\Gamma(A,A)$ and
$$f^{-1} f \equiv 1_{A \otimes \calk}$$
in $\Gamma(A \otimes \calk, A \otimes \calk)$.
That is, $f$ is invertible.

\end{itemize}

\subsection{Respecting split exactness}

\label{respecting_split_exactness}

The last characterization of $KK$-theory is split exactness.

\begin{itemize}

\item
For every split exact sequence (\ref{splitexact}) and $\vartheta_S$ as in \ref{uu} 
we have a diagram
\begin{equation}    \label{diagramvartheta}
\xymatrix{
& D \ar@<.5ex>[d]^{\vartheta_\cals} & \\
A \ar[ru]^{f} \ar@<.5ex>[r]^{i_A}
& A \oplus B \ar@<.5ex>[u]^{\sigma} \ar@<-.5ex>[r]_{p_B} \ar@<.5ex>[l]^{p_A} & B \ar@<-.5ex>[l]_{i_B}
\ar[lu]_{s}  ,
}
\end{equation}
where
\begin{equation}   \label{defsigma}
\sigma = p_A f + p_B s.
\end{equation}
For every such split exact sequence we define $\sigma$ to be invertible via $\vartheta_\cals$.
That is, we define the elementary equivalences 
$$\sigma \vartheta_\cals \equiv 1_{A \oplus B}$$ 
in $\Gamma(A\oplus B,A \oplus B)$ and
$$\vartheta_\cals \sigma \equiv 1_D$$
in $\Gamma(D,D)$.

\end{itemize}

\subsection{The operations respect equivalence}


We need to think about if the operations of taking the product and sums in $\Gamma$ is respected by
the equivalence relation $\equiv$ introduced in Section \ref{introducing_equivalences}.
That is, if $f_1 \equiv f_2$ and $g_1 \equiv g_2$ then we shall need to show that
$f_1 \times g_1 \equiv f_2 \times g_2$ and $\pm f_1 \pm g_1 \equiv \pm f_2 \pm g_2$.
Since the product in $\Lambda$ is bilinear by (\ref{defproduct}), and because of
the incorporation of the composed elementary equivalences (\ref{composedexpression}) to the elementary equivalences,
it is easy to check that all elementary equivalences we have introduced above satisfy this
for a given elementary equivalence $f_1 \equiv f_2$ and $g_1 = g_2$ (identity), or
if $f_1 = f_2$ (identity) and $g_1 \equiv g_2$ is
an elementary equivalence.
Hence the claim follows by using successive elementary equivalences.

\subsection{Definition of $GK$-theory}   

We are coming now to the definition of $GK$-theory. It is defined to be $\Gamma$ divided by equivalence.

\begin{definition}
{\rm
Let {\em generators and relations defined $KK$-theory $GK$} denote the category with object class $Obj(C^*)$ and morphism classes
$GK(A,B)$ to be
$\Gamma(A,B)$ divided by equivalence
defined in \ref{introducing_equivalences}-\ref{respecting_split_exactness}
for all objects $A$ and $B$ in $C^*$.
}
\end{definition}

\section{The universal property of $GK$-theory}

\label{section_universal_property}



$GK$-theory has the universal property described in this section. 

\begin{lemma}   \label{lemmaGKpropeties}
The category $GK$ is an additive category (if we accept that the morphism classes are classes and not sets).

The canonical functor
$$F:C^* \rightarrow GK : f \mapsto f$$
which maps objects 
identically to objects, and each morphism $f \in C^*(A,B)$ identically to the letter $f \in \Theta(A,B)$,
is split-exact, homotopy invariant and stable.

\end{lemma}

We shall discuss the last lemma. Of course each Hom-class 
$GK(A,B)$
is an abelian group (which is not small)
by \ref{introducing_equivalences}.
We have a zero object in $GK$ with the null $C^*$-algebra.
The reason is that the only existing morphism $1_{0 \oplus 0}$ in $C^*(0 \oplus 0,0 \oplus 0)$
must be the zero element of the Hom-class $GK(0 \oplus 0 , 0 \oplus 0)$ by identity (\ref{biproduct})
relative to the diagram
\begin{equation}
\xymatrix{0 \ar@<.5ex>[r]^{i_0} & 0 \oplus 0 \ar@<.5ex>[r]^{p_0} \ar@<.5ex>[l]^{p_0} & 0 \ar@<.5ex>[l]^{i_0} . }
\end{equation}
But the product with the zero element is always the zero element again,
so that $0 \oplus 0$ must be a zero object.
The product in $GK$ is bilinear by (\ref{defproduct}).
By (\ref{biproduct}) we have a biproduct and hence a coproduct in $GK$ by \cite[VIII.2.Thm. 2]{maclane}.
Hence $GK$ is an additive category.

The function
$F$ is a functor by \ref{respecting_homomorphisms} and \ref{the_unit_element}, which is homotopy invariant, stable and split-exact
by \ref{respecting_homotopy}, \ref{respecting_stability} and \ref{respecting_split_exactness}.

\begin{theorem}   \label{lemma_universal_GK}
Let $G: C^* \rightarrow E$ be a stable, homotopy invariant and split-exact functor from $C^*$ into an additive category $E$.
Then $G$ factorizes through $F$, that is, there exists a unique additive functor $\hat G$ such that the following diagram commutes:
\begin{equation}
\xymatrix{
& GK  \ar[rd]^{\hat G} & \\
C^* \ar[ru]^{F} \ar[rr]^{G} & & E .
}
\end{equation}
The functor $\hat G$ is given by
\begin{equation}  \label{ident5}
\hat G(A)=A
\end{equation}
for objects $A$ in $C^*$, 
%
\begin{equation}  \label{ident6}
\hat G(f) = G(f)
\end{equation}
for morphisms $f \in C^*(A,B)$, 
\begin{equation}   \label{ident1}
\hat G(f^{-1}) = G(f)^{-1}
\end{equation}
for a corner embedding $f:A \rightarrow A \otimes \calk$ and the letter $f^{-1}$ introduced
in \ref{uu}, and
\begin{equation}   \label{ident2}
\hat G(\vartheta_\cals) = \hat G(\sigma)^{-1} =  {\big (G(p_A) \circ G(f) + G(p_B) \circ G(s)\big )}^{-1}
\end{equation}
for every split exact sequence (\ref{splitexact}) and letter $\vartheta_\cals$ introduced in \ref{uu},
where $\sigma$ is defined in (\ref{defsigma}).

Moreover, one sets
\begin{equation}   \label{ident3}
\hat G(f_1 f_2 \ldots f_n) = G(f_1) \circ G(f_2) \circ \ldots \circ G(f_n)
\end{equation}
for $f_i$ in $\Theta(A_i,A_{i+1})$, and
\begin{equation}   \label{ident4}
\hat G(\pm f_1 \pm f_2 \pm \ldots  \pm f_n) = \pm G(f_1) \pm G(f_2) \pm \ldots \pm G(f_n)
\end{equation}
for $f_i$ in $\Lambda(A,B)$.

\end{theorem}

We are going to discuss this theorem. By (\ref{ident5}) and (\ref{ident6}) it is obvious that $\hat G \circ F = G$.
Because $G$ is a stable functor, for a corner embedding $f:A \rightarrow A \otimes \calk$ as in (\ref{ident1}) $G(f)$ is invertible. Also, every invertible element in a category
is uniquely defined. Hence identity (\ref{ident1}) is justified.

By the split-exactness of $G$, the split exact sequence (\ref{splitexactAB}) induces a coproduct diagram
\begin{equation}    \label{diagramxi}
\xymatrix{
 && G(A \oplus B) \ar@<.5ex>[dd]^{\xi^{-1}} \ar@<.5ex>[lldd]^{G(p_A)}  \ar@<-.5ex>[rrdd]^{G(i_B)}  & \\
  \\
G(A) \ar@<.5ex>[rruu]^{G(i_A)} \ar@<.5ex>[rr]^{i^A_E} && G(A) \sqcup G(B) \ar@<.5ex>[uu]^{\xi} \ar@<-.5ex>[rr]_{p^B_E} \ar@<.5ex>[ll]^{p^A_E} && G(B) \ar@<-.5ex>[ll]_{i^B_E} \ar@<-.5ex>[lluu]^{G(p_B)}  ,
}
\end{equation}
where $G(A) \sqcup G(B)$ denotes the coproduct and
$$\xi = p_E^A G(i_A) + p_E^B G(i_B)$$
is invertible. (We use the reversed order of notating composition of morphisms also in $E$.)
By \cite[VIII.2.Thm. 2]{maclane} the bottom line of (\ref{diagramxi}) is also a biproduct.
Notice that the diagram is commutative in the sense that
\begin{equation}   \label{identxi1}
i_E^A \xi = G(i_A)
\end{equation}
and
\begin{equation}   \label{identxi2}
\xi G(p_A) = p_E^A,
\end{equation}
and similarly so for $B$ instead of $A$.

By the split-exactness of $G$, the split exact sequence (\ref{splitexact}) induces a coproduct diagram
\begin{equation}    \label{diagramG}
\xymatrix{
& G(D) \ar@<.5ex>[d]^{\eta^{-1}}  & \\
G(A) \ar[ru]^{G(f)} \ar@<.5ex>[r]^{i^A_E}
& G(A) \sqcup G(B) \ar@<.5ex>[u]^{\eta} \ar@<-.5ex>[r]_{p^B_E} \ar@<.5ex>[l]^{p^A_E} & G(B) \ar@<-.5ex>[l]_{i^B_E}
\ar[lu]_{G(s)}  ,
}
\end{equation}
where
\begin{equation}   \label{defxi}
\eta = p_E^A G(f) + p_E^B G(s)
\end{equation}
is invertible. Entering here the identity 
(\ref{identxi2}) we get that
$$G(p_A) G(f) + G(p_B) G(s)$$
is invertible. This is, however, the element $\hat G(\sigma)$ occurring in (\ref{ident2}), and hence definition (\ref{ident2}) is valid.

Identity (\ref{ident3}) is necessary for a functor, 
and (\ref{ident4}) is necessary for a functor to be additive.
It is then clear by (\ref{defproduct}) that $\hat G$ respects products and sums.

In this way we at first obtain a well-defined (preliminary) function $\hat G: \Gamma \rightarrow E$.
We only need to check that equivalence in $\Gamma$ is respected by the function $\hat G$. That is, if $f_1 \equiv f_2$ then
$\hat G(f_1) \equiv \hat G(f_2)$. This is however true for all the elementary equivalences we have introduced, because
$G$ is a stable, homotopy invariant, split-exact functor into an additive category.
Only for equivalence (\ref{biproduct}) we need to remark that
an application of the functor $\hat G$ to the identity (\ref{biproduct}) and the use of identities (\ref{identxi1}) and (\ref{identxi2}) gives
\begin{equation}   \label{biproduct_E}
\xi^{-1} p^A_E i^A_E \xi +  \xi^{-1} p^B_E i^B_E \xi  =  1_{G(A \oplus B)}.
\end{equation}
This is, however, true because the bottom line of (\ref{diagramxi}) is a biproduct. 

\section{Turning Hom-classes to Hom-sets}

\label{section_homsets}

\subsection{Restricting the cardinality}

Our next aim is to turn the morphism classes into morphism sets.
To this end we restrict the cardinality of objects in $C^*$. That is, we choose a fixed cardinality $\chi$ and allow in $C^*$ only
those objects $A$ for which
$$card(A) \le \chi.$$

Then we select any {\em set} $\cala$ such that each object of $C^*$ is isomorphic to some object of $\cala$.
For example, we may choose the set of all $C^*$-subalgebras of $B(H)$ for some suitably large, fixed Hilbert space $H$.
We shall need a slightly larger object set, and enlarge $\cala$ to
the set
$$
\calb = \{ A, A \otimes \calk, \, A \oplus B \in Obj(C^*)\,|\, A,B \in \cala \}.$$


\subsection{First preparation with inverses}

Let us consider a corner embedding $f:A \rightarrow A \otimes \calk$ in $C^*$ and the letter $f^{-1}: A \otimes \calk \rightarrow A$ as introduced
in \ref{uu}. Choose an isomorphism $\pi: A \rightarrow B$ with $B \in \cala$.
Then we have a diagram
\begin{equation}
\xymatrix{
B \ar@<.5ex>[r]^{g} & B \otimes \calk  \ar@<.5ex>[l]^{g^{-1}}  
\\
A \ar[u]^\pi \ar@<.5ex>[r]^{f} & A \otimes \calk \ar[u]^{\pi \otimes id}  \ar@<.5ex>[l]^{f^{-1}} ,
}
\end{equation}
where $g$ is the unique morphism in $C^*$ such that the left square commutes, and $g^{-1}$ is the letter introduced in \ref{uu} for the corner embedding $g$.
As a path, $g$ is equivalent to $g \equiv \pi^{-1} f (\pi \otimes id)$ in $\Lambda(B,B \otimes \calk)$
(by (\ref{equivalencemorphism}) and $g= (\pi \otimes id) \circ f \circ \pi^{-1}$).
We take the inverse of $g$ and get
\begin{equation}    \label{substitutef}
f^{-1} \equiv (\pi \otimes id) g^{-1} \pi^{-1}.
\end{equation}
Notice that $g^{-1}$ is in $\Theta(B \otimes \calk,B)$ with source and range $B \otimes \calk,B \in \calb$. 


\subsection{Second preparation with inverses}

Similarly, given a split exact sequence (\ref{splitexact}) and $\vartheta$ as introduced in \ref{uu}, we consider
the diagram (\ref{diagramvartheta}) and double it by choosing isomorphisms $\pi_A:A \rightarrow A'$, $\pi_B:B \rightarrow B'$
and $\pi_D : D \rightarrow D'$ for which $A',B',D' \in \cala$.
We get a diagram like this:
\begin{equation}
\xymatrix{
 & & &   & D' \ar[d]^{\vartheta'} & \\
& D \ar[d]^{\vartheta} \ar@{.>}[rrru]^{\pi_D}  &  & A' \ar[ru]^{f'} \ar[r]^{i_A'} & A' \oplus B' \ar[u]^{\sigma'} \ar[r]^{p_B'} \ar[l]^{p_A'} & B' \ar[l]^{i_B'}
\ar[lu]_{s'}\\
A \ar[ru]^{f} \ar[r]^{i_A} \ar@{.>}[rrru]^{\pi_A} & A \oplus B \ar[u]^{\sigma} \ar[r]^{p_B} \ar[l]^{p_A}  \ar@{.>}[rrru]^{\pi_A \oplus \pi_B} & B \ar[l]^{i_B} \ar[lu]_{s}
\ar@{.>}[rrru]^{\pi_B}
}
\end{equation}

Here we define
\begin{eqnarray*}
&&\sigma = p_A f + p_B s,\\
&&\sigma' = p_A' f' + p_B' s' .
\end{eqnarray*}
The right hand sided triangle in the diagram and the letter $\vartheta'$ (from \ref{uu}) corresponds to the split exact sequence
\begin{equation}
\cals':\xymatrix{0 \ar[r] & A' \ar[r]^{f'} & D' \ar@<.5ex>[r]^{g'} & B' \ar[r] \ar@<.5ex>[l]^{s'} & 0  }
\end{equation}
taken from the split exact sequence (\ref{splitexact}) via the isomorphisms $\pi$.
We are going to show that the rectangle
spanned by the edges $\pi_D$ and $\pi_A \oplus \pi_B$ in the diagram commutes.
Indeed,
\begin{eqnarray*}
&&  \sigma \pi_D = (p_A f + p_B s) \pi_D = p_A f \pi_D + p_B s \pi_D
\equiv p_A \pi_A f' + p_B \pi_B s'  \\
&\equiv&  (\pi_A \oplus \pi_B) p_A' f' + (\pi_A \oplus \pi_B) p_B' s'
= (\pi_A \oplus \pi_B) \sigma'.
\end{eqnarray*}
Hence this rectangle involving the inverse $\vartheta$ instead of $\sigma$ commutes and we get
\begin{equation}    \label{substitutevartheta}
\vartheta \equiv  \pi_D \vartheta' (\pi_A \oplus \pi_B).
\end{equation}
Again we have achieved that $\vartheta'$ is in $\Theta(D',A' \oplus B')$ with
source and range
$D',A' \oplus B'$ lying in $\calb$.

\subsection{Rewriting words}

Now let us be given objects $A$ and $B$ in $C^*$ and an element in $\Lambda(A,B)$.
It is represented as a word $f_1 f_2 \ldots f_n$, or path (\ref{path}), with $f_i \in \Theta(A_i,A_{i+1})$.
In this word replace every letter $f_i$ of the form $f_i=f^{-1}$ or $f_i= \vartheta_\cals$ as introduced in \ref{uu} by the corresponding  
equivalent expression
(\ref{substitutef}) or (\ref{substitutevartheta}), respectively.
What comes out is a new word $g_1 g_2 \ldots g_m$ in $\Lambda(A,B)$ which is equivalent in $\Gamma(A,B)$ to the word $f_1 f_2 \ldots f_n$.
Notice that none of the synthetical letters $f^{-1}$ and $\vartheta_\cals$ as introduced in \ref{uu} follow each other in this new word, and each of these synthetical
letters has source and range in $\calb$. The letters between these synthetical letters are morphisms in $C^*$, and we fuse them together
by the equivalences (\ref{equivalencemorphism}).
The result is another word
$$h_1 q_1 h_2 q_2 \ldots h_{k-1} q_{k-1} h_k$$
which is equivalent to $g_1 g_2 \ldots g_m$, where each letter $h_i$ is a morphism in $C^*$ and each letter $q_i$ is a synthetical letter
$f^{-1}$ or $\vartheta_\cals$ as introduced in \ref{uu}.
Also, each $q_i$ is in $\Theta(D,D')$ for certain elements $D,D' \in \calb$.
Hence the path $f_1 f_2 \ldots f_n$ is equivalent to a path of the form
\begin{equation}   \label{pathA2}
\xymatrix{A \ar[r] & D_1 \ar[r] & D_3 \ar[r] &\ldots \ar[r] & D_{l} \ar[r] & B , }
\end{equation}
where $D_i \in \calb$.

\subsection{Rewriting sums}

If we have given an element of $\Gamma(A,B)$ then it is a formal sum (\ref{formalsum}) of elements $f_i$ in $\Lambda(A,B)$.
We may
apply the above procedure to each $f_i$ and so obtain the following lemma.

\begin{lemma}    \label{lemma_rewritting_morphisms}
Every element of $\Gamma(A,B)$ is equivalent to a formal sum
$$\pm f_1 \pm f_2 \pm \ldots \pm f_n$$
in $\Gamma(A,B)$ where
each $f_i$ is a path in $\Lambda(A,B)$ of the form (\ref{pathA2}) with $D_i \in \calb$.
\end{lemma}

\subsection{Equivalence with the category $\cald$}

Since each letter set $\Theta(A,B)$ is a set and $\calb$ is a set it is clear that the last lemma shows the following:

\begin{corollary}    \label{corollarydab}
For all objects $A$ and $B$ in $C^*$
there exists a set $\cald_{A,B}$ such that each element of $\Gamma(A,B)$ is equivalent to an element of $\cald_{A,B}$.
\end{corollary}

Say that two elements in $\cald_{A,B}$ are {locally-equivalent} if they are equivalent in $\Gamma(A,B)$.
It easy to see that this is an equivalence relation $R_{A,B}$. 
For all objects $A,B,C$ in $C^*$
define a product
$$\cald_{A,B}/R_{A,B} \; \times \; \cald_{B,C}/R_{B,C} \; \longrightarrow  \; \cald_{A,C}/R_{A,C}  \; : \;f \times g \mapsto h \equiv fg$$
by taking the product $f g$ in $\Gamma$ of two representatives $f \in \cald_{A,B}$ and $g \in \cald_{B,C}$ and then choosing
any $h \in \cald_{A,C}$ which is equivalent in $\Gamma$ to $f g$ (by Corollary \ref{corollarydab}).  
It is clear that this definition does not depend on the choices of the representatives $f,g$ and $h$.

Addition in $\cald_{A,B}/R_{A,B}$ we define in the same vein 
by using 
addition in $\Gamma(A,B)$.
Let us write
$$\cald(A,B) = \cald_{A,B}/R_{A,B}.$$
Let $\cald$ denote the category with objects $Obj(C^*)$ and these morphism sets.

This shows the following:

\begin{corollary}   \label{corollaryHomsets}
If we restrict the cardinality of the allowed objects $A$ in $Obj(C^*)$ to $card(A) \le \chi$ for some fixed $\chi$,
then the Hom-classes of $GK$, that is the classes $\Gamma(A,B)$ modulo equivalence, are actually Hom-sets.

More precisely, instead of the classes $\Gamma(A,B)$ and their notion of equivalence, product and sum we may equally well work with the
sets 
$\cald(A,B)$ and their notion of product and sum.  
\end{corollary}

More category-theoretically we 
may say:

\begin{corollary}
The categories $GK$ and $\cald$
are equivalent.
\end{corollary}

The functors of this equivalence are given by the identic embedding functor $\cald \rightarrow GK$ and the functor $GK \rightarrow \cald$
of rewriting morphisms as described in Lemma \ref{lemma_rewritting_morphisms}.  

\section{The isomorphism with $KK$-theory}  

\label{section_kk}

\subsection{Equivalence with Kasparov's $KK$-theory}

For separable $C^*$-algebras, $KK$-theory has the same universal property as $GK$ described in Theorem \ref{lemma_universal_GK} by N. Higson \cite{higson}.
Hence we get a commuting diagram
\begin{equation}
\xymatrix{
& GK  \ar@<.5ex>[rd]^{\hat G} & \\
C^* \ar[ru]^{F} \ar[rr]^{G} & & KK \ar@<.5ex>[lu]^{\hat F}.
}
\end{equation}
The functor $\hat F$ exists by Lemma \ref{lemmaGKpropeties}, Corollary \ref{corollaryHomsets} and \cite[Thm. 4.5]{higson}
for $C^*$-algebras without group action,
and \cite{thomsen} or \cite[Thm. 6.6]{meyer2000} for the inclusion of group actions. 
The aforementioned theorems also claim that both $\hat F$ and $\hat G$ are uniquely determined.

Since $\hat F \circ \hat G \circ F = \hat F \circ G = F$ and $\mbox{id}_{GK} \circ F = F$, the uniqueness assertion of
Theorem \ref{lemma_universal_GK} applied to $F$ shows that
$\hat F \circ \hat G = \mbox{id}_{GK}$.
Similarly we get $\hat G \circ \hat F = \mbox{id}_{KK}$.
Hence $\hat F$ and $\hat G$ are isomorphisms of categories and we get:

\begin{theorem}  \label{theoremKK}
If $Obj(C^*)$ is the class of all separable $M$-equivariant $C^*$-algebras 
then $GK$-theory 
is 
isomorphic with Kasparov's $KK$-theory. 
In other words, $GK \cong KK$.
\end{theorem}

\subsection{The descent homomorphism}

We shall remark how easy one may define a descent homomorphism (see \cite{kasparov1988}) by the universal property of $GK$-theory.
Let us now denote $M$-equivariant $KK$-theory by $KK^M$ and ordinary $KK$-theory (with trivial group action) by $KK$.
Similarly we write $C^*_M$ for $M$-equivariant $C^*$-category and $C^*$ for ordinary $C^*$-category.
It is well-known that the canonical functor
$$G:C^*_M \rightarrow KK : A \mapsto A \rtimes M$$
mapping a homomorphism $\varphi : A \rightarrow B$ to the morphism in $KK$ induced by the canonical homomorphism
$\varphi \rtimes \mbox{id} :A \rtimes M \rightarrow B \rtimes M$ is
stable, homotopy invariant and split exact.
%
Hence Theorem \ref{lemma_universal_GK} applied to $G$ and Theorem \ref{theoremKK} (for separable $C^*$-algebras)
immediately yield a {\em descent homomorphism}
$$\hat G: KK^M \rightarrow KK : A \mapsto A \rtimes M.$$
This works analogously for the reduced crossed product $A \rtimes_r G$ and yields an analogous
descent homomorphism
$$\hat G_r: KK^M \rightarrow KK : A \mapsto A \rtimes_r M.$$

\subsection{Other $KK$-theories}

As remarked in the introduction, the definition of $GK$-theory works also for other categories than $C^*$.
For example one could consider other algebras than $C^*$-algebras, or consider $C^*$-algebras
but equipped with an action by a semigroup, or a category, or an inverse semigroup and so on.
The homomorphisms should then be chosen to be equivariant in the respective sense. 
That the action is given by a group was only relevant in this section.

Differently, but closely related, N. Higson (see \cite{cuntzueberblick}),
and A. Connes and N. Higson \cite{conneshigson,conneshigson2} develop
a universal theory which is stable, homotopy invariant and half-exact. This theory is called $E$-theory
and the difference to the theory of this note is that here we have split-exactness instead of half-exactness.

In another direction J. Cuntz \cite{cuntz1997} develops bivariant $K$-theories especially for 
locally convex algebras.
The essential difference, beside the substitution of homotopies by diffotopies, to this paper is 
that in Cuntz' approach the theories are {half-exact} for short exact sequences with linear splits
and produce long exact sequences together with Bott periodicity. In this paper we have a {split exact} theory.

\bibliographystyle{plain}
\bibliography{references}

\end{document}